# THE DISTRIBUTION OF MAXIMA OF APPROXIMATELY GAUSSIAN RANDOM FIELDS

By Yuval Nardi,[1] David O. Siegmund[2,3] and Benjamin Yakir[3]

*Carnegie Mellon University, Stanford University and Hebrew University*

Motivated by the problem of testing for the existence of a signal of known parametric structure and unknown "location" (as explained below) against a noisy background, we obtain for the maximum of a centered, smooth random field an approximation for the tail of the distribution. For the motivating class of problems this gives approximately the significance level of the maximum score test. The method is based on an application of a likelihood-ratio-identity followed by approximations of local fields. Numerical examples illustrate the accuracy of the approximations.

**1. Introduction and summary.** There are two central themes in this paper. One is the development of a method for the derivation of analytic approximations for the tail of the distribution of the maximum of a smooth random field. Such random fields and the distribution of their maxima emerge naturally in a variety of statistical applications, for example, brain mapping or searching for hot spots of disease in space and/or time. See, for example, [10, 11, 12, 13, 14, 15, 16, 19].

The other theme involves the detailed investigation of a specific case, which is asymptotically Gaussian but where direct application of results for Gaussian fields does not seem adequate. This field arises in the context of testing for the presence of a signal of a given parametric structure within a noisy image. The image is composed of an array of pixels. The effect of a signal at a given pixel depends on the distance between the signal and the pixel. A score statistic is constructed for each candidate signal and an overall test statistic for the presence of some signal is obtained by maximizing the

Received May 2007; revised May 2007.
[1]Supported by a graduate grant from the Faculty of Social Sciences, Hebrew University.
[2]Supported by the National Science Foundation.
[3]Supported by a grant from the U.S.–Israel Binational Science Foundation (BSF 2000103).

score over the collection of all candidate "locations." We will consider second order approximations for the tail of the distribution of the test statistic under the null hypothesis of the absence of a signal—the significance level of the test. Similar methods can be applied to obtain an approximation for the power. In Section 5 we indicate how the method may also be adapted to other models, including (under different technical conditions) the frequently discussed case of smooth Gaussian fields, and how it can be used to obtain higher order approximations.

We begin by describing the random field of interest. Assume that there is a process $\{W_u : u \in \mathcal{A}_n\}$ over some space. These $|\mathcal{A}_n| = n$ observations are mutually independent with distributions $F_u$. The collection $\mathcal{A}_n$ indicates the locations of pixels, which may or may not be regularly spaced, and the index $u$ denotes the location of a particular pixel. The process of interest is $X_t = \sum_{u \in \mathcal{A}_n} \theta_u(t) W_u, t \in T$. We assume that $\theta : \mathcal{A}_n \times T \to \mathbb{R}^+$ is a known real-valued function and that $T$ is a nice subset of the $d$-dimensional Euclidean space. In typical applications $t$ is a parameterization of a putative signal and $\theta_u(t)$ represents the effect of $W_u$ on that signal. The case where $T$ is a rectangle with sides parallel to the coordinate axes, $t$ indexes subrectangles, also with sides parallel to the coordinate axes but having varying location and dimensions, and $\theta_u(t)$ is the indicator that $u$ belongs to $t$ was discussed by Siegmund and Yakir [14]. In this paper we consider the case where $\theta_u$ is a smooth function of $t \in T$ for each $u$. In an example that we consider later, $t$ is a line segment joining the left and right sides of the unit square and $\theta_u(t)$ is a decreasing function of the distance from the point $u$ to the line $t$. In this example, large values of $W_u$ for $u$ close to $t$ lead to large values of $X_t$ and indicate the presence of the signal $t$, while values of $W_u$ for $u$ far from $t$ have less effect on $X_t$.

We embed the distributions of $X_t, t \in T$ in an exponential family with a natural parameter $\xi$ by assuming that the $W_u$'s obey an exponential distribution laws of the form:

$$(1.1) \quad dF_{t,\xi}(w_u) = \exp\{\xi\theta_u(t)w_u - \psi_u(\xi\theta_u(t))\} dF_u(w_u), \qquad \xi \in \mathbb{R}, t \in T.$$

Throughout the paper we use the convention that derivatives with respect to $\xi$ are denoted by apostrophes, while derivatives with respect to $t$ by dots. We assume that the distributions have been standardized so that

$$(1.2) \qquad \psi_u(0) = 0, \qquad \psi'_u(0) = 0, \qquad \psi''_u(0) = 1.$$

We can formulate the null hypothesis of no signal as $H_0 : \xi = 0$ for all $t$, while under the alternative $\xi = \xi_t > 0$ for some value of $t$. In terms of the exponential embedding, the log-likelihood for fixed $t$ is given by

$$(1.3) \qquad l_n(t, \xi) = \sum_{u \in \mathcal{A}_n} \{\xi\theta_u(t)W_u - \psi_u(\xi\theta_u(t))\}.$$

Differentiating the log-likelihood (1.3) twice with respect to $\xi$ and substituting $\xi = 0$, we obtain the standardized score statistic at a fixed $t \in T$, which by virtue of (1.2) is given by

$$(1.4) \qquad Z_n(t) = \frac{l'_n(t,0)}{\sqrt{I_n(t)}} = \sum_{u \in \mathcal{A}_n} \beta_{u,n}(t) W_u,$$

where $I_n(t) = \sum_{u \in \mathcal{A}_n} \theta_u^2(t)$ is the Fisher information and $\beta_{u,n}(t) = \theta_u(t)/[I_n(t)]^{1/2}$. A (one-sided) test statistic is obtained by maximizing $Z_n(t)$ over $T$. Since the value of $t$ giving rise to the signal is unknown, we consider as a test statistic $\max_t Z_n(t)$. The associated $p$-value is given by the probability:

$$(1.5) \qquad \mathbb{P}\left(\sup_{t \in T} Z_n(t) \geq x\right)$$

computed under the assumption that $\xi = 0$.

Consider the random field $\{Z_n(t), t \in T\}$. For each fixed $t$, when $\xi = 0$ the random variable $Z_n(t)$ is asymptotically standard normal. The sample points of this field are real valued functions on $T$, which are smooth functions of $t$, since the $\{\theta_u(t)\}$ are. The main result that we would like to establish is that when $x = o(n^{1/4})$, the probability (1.5) can be approximated, up to a term inside the braces that is $o(1/x)$ by

$$(1.6) \quad \begin{aligned} &x^{d-1}(2\pi)^{-d/2}\phi(x) \\ &\quad \times \bigg\{ \int_T e^{-\delta_n(t)} |\Lambda_n(t)|^{1/2} (1 - r_n^2(t)/2\sigma_n^2(t))\, dt \\ &\quad + \frac{1}{x}\left(\frac{\pi}{2}\right)^{1/2} \\ &\quad \times \int_{\partial T} e^{-\delta_n(t)} (|\Lambda_n(t)| \cdot \langle \dot{g}(t), \Lambda_n^{-1}(t)\dot{g}(t)\rangle)^{1/2}/\|\dot{g}(t)\|\, dV_{\partial T}(t) \bigg\}. \end{aligned}$$

The indicated approximation involves powers of the threshold $x$ and the standard normal density $\phi$. It also involves integration of functions of $t$ denoted by $\Lambda_n(t)$, $\delta_n(t)$, $r_n(t)$, $\sigma_n^2(t)$ and $\dot{g}(t)$, which are defined below. Angular brackets "$\langle \cdot, \cdot \rangle$" correspond to the inner product of vectors and "$\|\cdot\|$" to the Euclidean norm.

The region $T$ is defined with the aid of constraint functions:

$$(1.7) \qquad T = \{t \in \mathbb{R}^d : g_i(t) \leq 0, 1 \leq i \leq m\}.$$

These functions are generically denoted $g(t)$, with gradient vector $\dot{g}(t)$. The boundary of the region, $\partial T$, corresponds to values of the parameters, in general a manifold, where $g(t) = 0$. The differential $dV_{\partial T}$ denotes the volume element of the ($d - 1$ dimensional part of the) boundary of $T$.

The other functions in (1.6) are functionals of $Z_n(t)$, the value of the field at $t$, and of $\dot{Z}_n(t)$, the random gradient of the field at that point. These functionals are computed under the alternative distribution for $W_u$ given by (1.1) with the amplitude

$$\xi_t = \xi_{t,n,x} = x I_n^{-1/2}(t). \tag{1.8}$$

Specifically, $\Lambda_n(t)$ is a positive-definite matrix defined by

$$\Lambda_n(t) = \frac{\sum_{u \in \mathcal{A}_n}[\dot{\theta}_u(t) \otimes \dot{\theta}_u(t)]}{I_n(t)} - \frac{[\dot{I}_n(t) \otimes \dot{I}_n(t)]}{4 I_n^2(t)}, \tag{1.9}$$

where "$\otimes$" is the outer (Kronecker) product. The matrix $\Lambda_n(t)$ is (asymptotically) the covariance matrix of $\dot{Z}_n(t)$. It is related to the expected value of the Hessian of $Z_n(t)$ under the tilted measure by $\sum_{u \in \mathcal{A}_n} \ddot{\beta}_{u,n}(t) \psi'(\xi_t \theta_u(t)) = -x \Lambda_n(t) + O(x^2 n^{-1/2})$.

The term $\delta_n(t)$ is the difference between the actual log-likelihood and the log-likelihood of the approximating normal distributions:

$$\delta_n(t) = l_n(t, \xi_t) - [x Z_n(t) - x^2/2] = x^2/2 - \sum_{u \in \mathcal{A}_n} \psi_u(\xi_t \theta_u(t)). \tag{1.10}$$

This quantity is deterministic and is of the order of magnitude $O(x^3 n^{-1/2})$. Another measure of discrepancy from the normal limit is $r_n(t)$, which is the difference between the (tilted) expectation of $Z_n(t)$ and the threshold $x$:

$$r_n(t) = \sum_{u \in \mathcal{A}_n} \beta_{u,n}(t) \psi_u'(\xi_t \theta_u(t)) - x. \tag{1.11}$$

This term is of the order of magnitude of $O(x^2 n^{-1/2})$. Finally, the term $\sigma_n^2(t)$ is given in terms of the $d \times d$ covariance matrix of the gradient:

$$\Sigma_n(t) = \sum_{u \in \mathcal{A}_n} [\dot{\beta}_{u,n}(t) \otimes \dot{\beta}_{u,n}(t)] \psi_u''(\xi_t \theta_u(t))$$

and the correlation between $Z_n(t)$ and $\dot{Z}_n(t)$:

$$\rho_n(t) = \sum_{u \in \mathcal{A}_n} [\beta_{u,n}(t) \dot{\beta}_{u,n}(t)] \psi_u''(\xi_t \theta_u(t))$$

by

$$\sigma_n^2(t) = 1 - \langle \rho_n(t), \Sigma_n^{-1}(t) \rho_n(t) \rangle. \tag{1.12}$$

Observe that in the Gaussian limit $Z(t)$ and $\dot{Z}(t)$ become independent, and hence $\sigma_n(t)$ converges to one.

Sufficient regularity conditions for (1.6) are stated in Theorem 4.9 below. In Section 2 we outline the principles of the method for producing expansions of the probability (1.5) in smooth random fields. Section 3 presents some numerical examples. In Section 4 we provide details regarding the approximation of the significance level, with some proofs deferred to an Appendix. Section 5 discusses some extensions.

REMARKS. (i) There is a large related literature in the case of Gaussian fields, for example, [3, 4, 7, 13, 16, 17]. See [1] for a review up to 2000 and additional references, and [18] for an outstanding recent contribution. To our knowledge the non-Gaussian case is relatively unexplored, and even then the published results are for relatively simple concrete problems, and are derived heuristically. See, for example, [11] and [10].

(ii) In the first integral in (1.6), the differential $|\Lambda_n(t)|^{1/2}\,dt$ is easily recognized as the volume element for a manifold with metric tensor $\Lambda_n$, exactly as one knows it must be from the familiar case of a Gaussian random field. Similarly, one knows from the Gaussian case that the appropriate differential in the second integral (i.e., the product of all factors in the integrand *except* the factor $e^{-\delta_n(t)}$) must equal the volume element for the boundary of the manifold with metric tensor $\Lambda_n$. It seems easy to prove this result in special cases, although difficult in general. [E.g., one can show the equivalence when the boundary of $T$ is given (locally) by sets like the set where $g(t) = t_d - f(t_1,\ldots,t_{d-1})$ equals 0, for a suitable smooth function $f$.] However, since application of equation (1.6) does not rely on this interpretation, we do not pursue the argument here.

(iii) In exponential change of measure arguments, one often chooses the parameter, here $\xi$, so that under the new distribution the process of interest has expectation exactly equal to the threshold $x$. Since this would result in a nonlinear equation for $\xi$, we find it simplifies some Taylor series expansions to use a slightly different value, which has an explicit form, and would be the conventional value if the process were exactly Gaussian.

**2. Outlining the method.** Our method for expanding probabilities of the form (1.5) is based on a application of a likelihood ratio identity, followed by local expansions. It is very similar to the approach that has been applied in previous work, which was, however, related to fields of Brownian-motion type. See [14, 15], and [20]. However, in the current work some steps have been modified in order to exploit the smoothness of the random paths. The application of the methodology is split into six building blocks, which are described below. Details are given in Section 4.

**Measure transformation:** The first step involves a likelihood ratio identity. This transformation allows us to recenter the analysis in a setting where the probability of crossing the threshold is substantially larger and where the central limit theorem is more likely to be applicable.

**Localization:** The likelihood ratio identity produces a functional of the random field. Conditioning on the signal, that is, the parameter point selected by the alternative distribution, it is argued that the value of the functional is determined mainly by the local behavior of the field about that point in the parameter space, so the original field can be replaced by a finite order Taylor polynomial and the computation of the functional restricted to a smaller subset of the parameter space.

**Application of a central limit theorem:** The local term that emerges from the previous step is a functional of the value of the process at the signal location and of derivatives of the field at that point. The functional is approximated by a similar functional applied to an appropriate multinormal distribution.

**Elimination of the indicator:** The functional has an exponential component, a component in the form of an indicator and components that are essentially polynomial in the derivatives of the random field. By conditioning on these derivatives and using a Mill's ratio type of approximation, one can eliminate the indicator.

**Evaluation of the functional:** In theory this step, which involves approximation of the functional of the derivatives by polynomials and the evaluation of the expectation of the resulting Gaussian polynomials, is straightforward. In practice, it is tedious to apply if higher order approximations are sought. It is at this stage that boundary effects become significant.

**Integration over the parameter space:** The analysis in the last four steps is carried out at the signal location. A final integration of the functionals over the parameter space must be carried out to produce an approximation to the probability of interest.

**3. Simulation studies.** In this section we examine via simulations the accuracy of (1.6) for one example of an approximately Gaussian random field. The simulations were programmed using the C++ language. Probabilities were approximated using 5,000 iterations of the simulations, which corresponds to results accurate up to the second digit after the decimal point.

Consider a background field of independent Bernoulli random variables covering the standard unit square on a fixed dyadic grid $\mathcal{A}_n = \{(i/2^m, j/2^m), 0 \leq i, j \leq 2^m\}$, so $n = (2^m + 1)^2$. A signal is composed of a straight line that passes between a point on the intersection of the unit square with the vertical line $x = 0$ and a point on its intersection with the vertical line $x = 1$. This figure can be imagined to be a very primitive "edge" in a two dimensional image. More interesting examples along these lines would involve, say, broken lines with some maximum number of breaks at unspecified positions. The parameter space $T$ is again the unit square, with each parameter point representing the vertical levels of the leftmost and the rightmost points of intersection, respectively. For the function $\theta_u(t)$, which measures the closeness of a point $u \in \mathcal{A}_n$ to the signal $t \in T$, we use $\theta_u(t) = \exp[-D \cdot d(u,t)^2/2]$, where $d(u,t) = \min_{0 \leq p \leq 1} \|u - [p(0, t_1) + (1-p)(1, t_2)]\|$ and $D$ is a positive scale parameter. The random variables $W_u$ were taken to be Bernoulli with probability of success $p_0 = 0.1$. When a signal is present, this probability is shifted to $p_1 = p_0/[p_0 + (1 - p_0)e^{-\xi \theta_u(t)}]$.

TABLE 1
*Comparison between a Gaussian approximation ($p_G$) and the approximation ($p_E$) based on (1.6). Empirical significance level ($\hat{p}$) was estimated from 5,000 iterations (SD $\in [0.0015, 0.003]$)*

| | $D = 10$ | | | | $D = 20$ | | | | $D = 50$ | | |
|---|---|---|---|---|---|---|---|---|---|---|---|
| $x$ | $\hat{p}$ | $p_E$ | $p_G$ | $x$ | $\hat{p}$ | $p_E$ | $p_G$ | $x$ | $\hat{p}$ | $p_E$ | $p_G$ |
| 2.5 | 0.046 | 0.036 | 0.026 | 2.8 | 0.044 | 0.043 | 0.023 | 3.1 | 0.042 | 0.088 | 0.024 |
| 2.6 | 0.036 | 0.029 | 0.020 | 2.9 | 0.034 | 0.035 | 0.018 | 3.2 | 0.036 | 0.075 | 0.018 |
| 2.7 | 0.029 | 0.024 | 0.016 | 3.0 | 0.030 | 0.029 | 0.014 | 3.3 | 0.034 | 0.064 | 0.013 |
| 2.8 | 0.021 | 0.020 | 0.012 | 3.1 | 0.021 | 0.024 | 0.010 | 3.4 | 0.024 | 0.054 | 0.010 |
| 2.9 | 0.015 | 0.016 | 0.009 | 3.2 | 0.015 | 0.019 | 0.008 | 3.5 | 0.020 | 0.046 | 0.007 |
| 3.0 | 0.013 | 0.013 | 0.007 | 3.3 | 0.014 | 0.016 | 0.006 | 3.6 | 0.013 | 0.040 | 0.005 |
| 3.1 | 0.009 | 0.010 | 0.005 | 3.4 | 0.012 | 0.013 | 0.004 | 3.7 | 0.010 | 0.034 | 0.004 |

Table 1 compares, with $m = 5$ ($n = 1{,}089$), the simulated tail probability ($\hat{p}$) with the analytic approximation ($p_E$) in (1.6), and a parallel approximation based on treating the field as Gaussian ($p_G$), for which (1.6) is used, but with $\delta_n = 0 = r_n$. Three values, $D = 10, 20, 50$ and a range of thresholds $x$ corresponding to $p$-values that vary between 0.05 and 0.01 were examined. Overall, when the values of the scale parameter are not too large ($D = 10, 20$) the approximations given by (1.6) produce good results and are more accurate than an approximation based on the Gaussian limit. When $D = 50$, only a small fraction of the background observations contribute to the score statistic, so the marginal distribution of $Z(t)$ is poorly approximated by a normal distribution. The Gaussian approximation is anticonservative, while (1.6) overcompensates and is very conservative. It would be interesting to know whether the large deviation like approximation suggested in [10] can be adapted to deal with this case.

In Table 2 we compare the approximated $p$-values for $m = 5$ ($n = 1{,}089$) versus $m = 6$ ($n = 4{,}225$) when $D = 17$. Since the accuracy of the approximations is comparable in both cases one is tempted to conclude that (1.6) is stable with respect to the sample size in the range where the approximation is valid.

**4. Detailed proofs.** In this section we add the details to the outline that was presented in Section 2. Throughout, we will try to keep the discussion as general as possible in order to lay the foundation for the extension of the method to other models and to higher order approximations.

4.1. *Measure transformation.* We begin by transforming the null probability measure under which the field is centered and the probability in (1.5) is relatively small to an alternative probability.

Put a uniform prior over the parameter space $T$. Let $\mathbb{P}_t \equiv \mathbb{P}_{t,\xi_t}$ correspond to the probability measure for which, given $t$, the random variables $\{W_u\}$ are (conditionally) independent and are distributed according to (1.1) with parameter $\xi_t \theta_u(t)$, where $\xi_t = xI_n^{-1/2}(t)$ was defined in (1.8). Let $\mathbb{P}_x$ be the unconditional distribution of $\{W_u\}$, so

$$(4.1) \qquad \mathbb{P}_x(\cdot) = \frac{1}{\lambda(T)} \int_T \mathbb{P}_t(\cdot)\, dt,$$

for $\lambda(T)$ the Lebesgue measure of $T$.

The log-likelihood ratio of $\mathbb{P}_t$ relative to $\mathbb{P}$ is given by (1.3) with $\xi_t$ in place of $\xi$. Hence, the likelihood ratio of $\mathbb{P}_x$ relative to $\mathbb{P}$ can be written as

$$(4.2) \qquad \frac{d\mathbb{P}_x}{d\mathbb{P}} = \frac{1}{\lambda(T)} \int_T e^{l_n(s,\xi_s)}\, ds.$$

Applying a likelihood ratio identity, one obtains

$$\mathbb{P}\left(\sup_{t \in T} Z_n(t) \geq x\right) = \int_T \mathbb{E}_t\left[\frac{1}{\int_T e^{l_n(s,\xi_s)}\, ds}; \sup_{s \in T} Z_n(s) \geq x\right] dt.$$

Simple algebraic manipulations lead to the representation in the form

$$(4.3) \quad = x^d e^{-x^2/2} \int_T e^{-\delta_n(t)} \mathbb{E}_t\left[\frac{e^{-x(Z_n(t)-x)}}{S}; x(Z_n(t)-x) + \log M \geq 0\right] dt,$$

where for each $t$, $M$ and $S$ are defined in terms of the random field $X_t(s) = X_{t,x,n}(s) = x(Z_n(s) - Z_n(t)), s \in T$, by the equations $M = \sup_{s \in T} \exp[X_t(s)]$ and $S = \int_{s \in T} x^d \exp[X_t(s) + \delta_n(s) - \delta_n(t)]\, ds$, with $\delta_n(\cdot)$ as defined in (1.10).

TABLE 2
*Comparison of the approximation ($p_E$) between $m = 5$ and $m = 6$ when $D = 17$. Empirical significance level ($\hat{p}$) was estimated from 5,000 iterations (SD $\in [0.0015, 0.003]$)*

|       | $m = 5$ |       | $m = 6$ |       |
|-------|---------|-------|---------|-------|
| $x$   | $\hat{p}$ | $p_E$ | $\hat{p}$ | $p_E$ |
| 2.6   | 0.049   | 0.054 | 0.045   | 0.043 |
| 2.7   | 0.040   | 0.044 | 0.037   | 0.034 |
| 2.8   | 0.032   | 0.036 | 0.031   | 0.027 |
| 2.9   | 0.029   | 0.030 | 0.020   | 0.022 |
| 3.0   | 0.020   | 0.024 | 0.019   | 0.017 |
| 3.1   | 0.019   | 0.020 | 0.015   | 0.013 |
| 3.2   | 0.013   | 0.016 | 0.009   | 0.010 |
| 3.3   | 0.009   | 0.013 | 0.006   | 0.008 |

4.2. *Localization and a Taylor approximation.* In this section we concentrate on the integrand in representation (4.3). The term $M$ produced by maximization and the term $S$ produced by integration will be replaced by similar, but more tractable, expressions. The replacement is two-fold. First, a subset of the parameter space is used instead of the complete parameter space. This subset, which we denote by $V_t$, is a neighborhood of $t$ with radius depending on the threshold $x$ and on the sample size $n$. The resulting local random field is denoted by $\mathcal{X}_t = \{X_t(s), s \in V_t\}$. Second, this local field is replaced by its Taylor approximation about $t$ (to a certain order). Denote $\hat{\mathcal{X}}_t = \{\hat{X}_t(s), s \in V_t\}$ to be the random field generated by these two approximations. The main subject in this subsection is to show that by taking the radius of $V_t$ to be of an appropriate order and by taking enough terms in the Taylor expansion of the field, one may control the error involved up to a pre-prescribed level of accuracy.

For $r > 0$, let $B(r)$ be the $d$-ball of radius $r$, as defined by the Euclidean metric. Denote the elements of any $d$-dimensional vector by subscripts, for example, $s = (s_1, \ldots, s_d)$. For every $k$-tuple $i_1, \ldots, i_k$ from $\{1, \ldots, d\}$, and a $d$-dimensional vector $x$ let $x_{i_1,\ldots,i_k} = \prod_{j=1}^{k} x_{i_j}$. For every function on $T$ and on $T \times T$ we denote partial derivatives using superscripts. For example, $\beta_{u,n}^{ij}(t)$ is the second order partial derivative of $\beta_{u,n}(s)$ with respect to $s_i, s_j$, evaluated at $s = t$. Throughout, we use the Einstein summation convention (cf. [5], pages 136–137).

We will consider the case where the information regarding $\xi$ at $t$, $I_n(t) = \sum_{u \in \mathcal{A}_n} \theta_u^2(t)$, is of the order of magnitude of the sample size $n$. This assumption clearly holds for a functions $\theta$ which are uniformly bounded and bounded away from zero over $\mathcal{A}_n \times T$. We assume henceforth that uniform boundedness holds both for $\theta$ and for its derivatives, up to an appropriate order. We also assume boundedness of $\psi$ and its derivatives.

By (1.4) we can decompose $Z_n(s)$ as

$$(4.4) \qquad Z_n(s) = U_n(s) + \sum_{u \in \mathcal{A}_n} \beta_{u,n}(s) \psi_u'(\xi_t \theta_u(t)),$$

where $U_n(s) = \sum_{u \in \mathcal{A}_n} \beta_{u,n}(s)(W_u - \psi_u'(\xi_t \theta_u(t)))$. Observe that $U_n(s)$ is a centered variable under the probability measure $\mathbb{P}_t$. The second component is deterministic. The investigation of the localized field $X_t(s) = x(Z_n(s) - Z_n(t))$ amounts to expanding $\beta_{u,n}(s)$ in a neighborhood of $t$. One needs to evaluate the derivatives of $\beta_{u,n}(t)$ explicitly in order to proceed with the analysis. The derivatives are incorporated separately with respect to the random part and the deterministic part.

Consider in detail the first two derivatives in conjunction with the deterministic part. Other derivatives (random and deterministic) will have, by the choice of the neighborhood, less of an effect on the overall value of $X_t(s)$. Direct computation gives

$$\dot{\beta}_{u,n}(s) = \frac{\dot{\theta}_u(s)}{I_n^{1/2}(s)} - \frac{J_n(s)}{I_n^{3/2}(s)}\theta_u(s),$$

$$\ddot{\beta}_{u,n}(s) = \frac{\ddot{\theta}_u(s)}{I_n^{1/2}(s)} - 2\frac{\dot{\theta}_u(s) \otimes J_n(s)}{I_n^{3/2}(s)}$$

$$- \frac{\theta_u(s)\dot{J}_n(s)}{I_n^{3/2}(s)} + 3\frac{\theta_u(s)J_n(s) \otimes J_n(s)}{I_n^{5/2}(s)},$$

where $J_n(s) = \sum_{u \in \mathcal{A}_n} \theta_u(s)\dot{\theta}_u(s)$ [so that $\dot{I}_n(s) = 2J_n(s)$].

In the following we let $\Theta(y)$ denote a function that is bounded above and below by expressions of the form $\text{const} \times y$.

By assumption, the first derivative term $\sum_{u \in \mathcal{A}_n} \dot{\beta}_{u,n}(t)\psi'(\xi_t\theta_u(t))$ is of the order of magnitude of $O(x^2 n^{-1/2})$, and is denoted by $r_1$. The second derivative term, $\sum_{u \in \mathcal{A}_n} \ddot{\beta}_{u,n}(t)\psi'(\xi_t\theta_u(t))$, can be written as $-x\Lambda_n(t) + r_2$, where $r_2$ is an $O(x^2 n^{-1/2})$ term, and $\Lambda_n(t)$ is a $\Theta(1)$, nonnegative definite matrix, given in (1.9).

Let $\alpha \geq 0$ and set

(4.5) $$V_t = t \oplus \Lambda_n^{-1/2}(t)B(\log x/x).$$

Each element of the approximating field $\{\hat{X}_t(s), s \in V_t\}$ is defined via a finite expansion:

$$\hat{X}_t(s) = x\langle s - t, \dot{U}_n(t)\rangle - \frac{x^2}{2}\langle s - t, \Lambda_n(t)(s - t)\rangle + r_1 + r_2$$

(4.6)
$$+ \sum_{k=2}^{\alpha+1} \frac{1}{k!} x(s-t)_{i_1,\ldots,i_k} U_n^{i_1,\ldots,i_k}(t)$$

$$+ \sum_{k=3}^{\alpha+2} \frac{1}{k!} x^2(s-t)_{i_1,\ldots,i_k} \sum_{u \in \mathcal{A}_n} \beta_{u,n}^{i_1,\ldots,i_k}(t)\psi'_u(\xi_t\theta_u(t))/x.$$

In the case $\alpha = 0$, only the first line of (4.6) is used. For $\alpha = 1$, which is relevant for the results presented in Section 1, one term from the second line and one term from the third line of (4.6) are also required.

Denote the remainder of the Taylor expansion by

(4.7) $$r_t(s) = r_{t,x,n}(s) = X_t(s) - \hat{X}_t(s).$$

Define

$$S_0 = \int_{s \in V_t} x^d e^{X_t(s) + \delta_n(s) - \delta_n(t)} \, ds, \qquad M_0 = \sup_{s \in V_t} e^{X_t(s)},$$

$$\hat{S}_0 = \int_{s \in V_t} x^d e^{\hat{X}_t(s) + \delta_n(s) - \delta_n(t)} \, ds, \qquad \hat{M}_0 = \sup_{s \in V_t} e^{\hat{X}_t(s)}.$$

We also set $X = x(Z_n(t) - x)$. The main result of this subsection is:

THEOREM 4.1. *Assume that $T$ is compact. Suppose further that $\theta_u(\cdot)$ belongs to $\mathbb{C}^{\alpha+3}$, and that all derivatives up to this order are bounded from above in $V_t$. Assume also that $\psi_u(\cdot)$ is three times differentiable with bounded derivatives over an open interval that contains the origin. Finally, assume that $I_n(t)$ is $\Theta(n)$, the smallest eigenvalue of $\Lambda_n(t)$ is bounded away from zero, and the largest eigenvalue is finite. Take $\varepsilon = \Theta(x^{-\alpha})$. Then*

$$\mathbb{E}_t\left[\frac{e^{-X}}{S}; X + \log M \geq 0\right] \leq e^{\varepsilon} \mathbb{E}_t\left[\frac{e^{-X}}{\hat{S}_0}; X + \log \hat{M}_0 \geq -\varepsilon\right] + o(x^{-\alpha}),$$

$$\mathbb{E}_t\left[\frac{e^{-X}}{S}; X + \log M \geq 0\right] \geq \frac{e^{-\varepsilon}}{1+\varepsilon} \mathbb{E}_t\left[\frac{e^{-X}}{\hat{S}_0}; X + \log \hat{M}_0 \geq \varepsilon\right] + o(x^{-\alpha}).$$

REMARKS. In the proofs given in this subsection we do not mention the contribution of $\delta_n(t)$ explicitly. In many interesting cases, which include for example the classical medium deviation of a normalized sum of i.i.d. variables, one considers threshold levels $x$ that satisfy $x = o(n^{1/6})$. For such $x$, $\delta_n(s)$ is a remainder term which tends to zero as $x$ and $n$ tend to infinity. Therefore, ignoring it throughout the various lemmas below will present no loss of generality. In other applications, where $x$ may tend to infinity slightly faster, for example, $x = o(n^{1/4})$, an appropriate choice of the constant $C$ in Lemma 4.2 and the bound over $M_0/S_0$, given in (4.17) below allow one to dominate $\delta_n(s)$ and proceed with the proof. We reintroduce $\delta_n(s) - \delta_n(t)$ in Section 4.5, where its exact contribution becomes relevant.

PROOF OF THEOREM 4.1. Roughly speaking, the theorem states that the error committed by switching between $S, M$ and $\hat{S}_0, \hat{M}_0$, respectively, is sufficiently small for our subsequent calculations. Therefore, one must make sure that (i) the remainder is small, and (ii) the contribution of quantities outside $V_t$ are negligible. To this end, the region $T \setminus V_t$ is covered with $K$ balls of radius $1/x^2$. Let $V_i$, for $1 \leq i \leq K$, denote a generic ball centered at $\tau_i \in T \setminus V_t$. Let $S_i$ and $M_i$ denote the analogues of $S_0$ and $M_0$ for such a ball. Note that $M = \max\{M_0, M_1, \ldots, M_K\}$ and $S \geq \max\{S_0, S_1, \ldots, S_K\}$. Since $M_0 \geq 1$, we see that $M_i \geq 1$ on $\{M = M_i\}$. The remainder of the proof consists of obtaining suitable upper and lower bounds.

The upper bound is obtained by a localization argument and a Taylor approximation (in that order):

$$\mathbb{E}_t\left[\frac{e^{-X}}{S}; X + \log M \geq 0\right]$$

$$\leq \mathbb{E}_t\left[\frac{e^{-X}}{S}; X + \log M \geq 0, M = M_0\right]$$

$$+ \sum_{i=1}^{K} \mathbb{E}_t\left[\frac{e^{-X}}{S}; X + \log M \geq 0, M = M_i\right]$$

$$\leq \mathbb{E}_t\left[\frac{e^{-X}}{S_0}; X + \log M_0 \geq 0\right] + \sum_{i=1}^{K} \mathbb{E}_t\left[\frac{M_i}{S_i}; M_i \geq 1\right]$$

$$= \mathbb{E}_t\left[\frac{e^{-X}}{S_0}; X + \log M_0 \geq 0, \sup_{s \in V_t} |r_t(s)| \leq \varepsilon\right]$$

$$+ \mathbb{E}_t\left[\frac{e^{-X}}{S_0}; X + \log M_0 \geq 0, \sup_{s \in V_t} |r_t(s)| > \varepsilon\right] + \sum_{i=1}^{K} \mathbb{E}\left[\frac{M_i}{S_i}; M_i \geq 1\right]$$

$$\leq e^\varepsilon \mathbb{E}_t\left[\frac{e^{-X}}{\hat{S}_0}; X + \log \hat{M}_0 \geq -\varepsilon\right] + \mathbb{E}_t\left[\frac{M_0}{S_0}; \sup_{s \in V_t} |r_t(s)| > \varepsilon\right]$$

$$+ \sum_{i=1}^{K} \mathbb{E}_t\left[\frac{M_i}{S_i}; M_i \geq 1\right].$$

The lower bound is obtained by going the other way:

$$\mathbb{E}_t\left[\frac{e^{-X}}{\hat{S}_0}; X + \log \hat{M}_0 \geq \varepsilon\right]$$

$$= \mathbb{E}_t\left[\frac{e^{-X}}{\hat{S}_0}; X + \log \hat{M}_0 \geq \varepsilon \geq \sup_{s \in V_t} |r_t(s)|\right]$$

$$+ \mathbb{E}_t\left[\frac{e^{-X}}{\hat{S}_0}; X + \log \hat{M}_0 \geq \varepsilon, \sup_{s \in V_t} |r_t(s)| > \varepsilon\right]$$

$$\leq e^\varepsilon \mathbb{E}_t\left[\frac{e^{-X}}{S_0}; X + \log \hat{M}_0 \geq \varepsilon, \sup_{s \in V_t} |r_t(s)| \leq \varepsilon\right]$$

$$+ e^{-\varepsilon} \mathbb{E}_t\left[\frac{\hat{M}_0}{\hat{S}_0}; \sup_{s \in V_t} |r_t(s)| > \varepsilon\right]$$

$$\leq e^\varepsilon \mathbb{E}_t\left[\frac{e^{-X}}{S_0}; X + \log M_0 \geq 0\right] + e^{-\varepsilon} \mathbb{E}_t\left[\frac{\hat{M}_0}{\hat{S}_0}; \sup_{s \in V_t} |r_t(s)| > \varepsilon\right]$$

$$= e^\varepsilon \Bigg( \mathbb{E}_t\left[\frac{e^{-X}}{S_0}; X + \log M_0 \geq 0, S \leq (1+\varepsilon)S_0\right]$$

$$+ \mathbb{E}_t\left[\frac{M_0}{S_0}; S > (1+\varepsilon)S_0\right]\Bigg)$$

$$+ e^{-\varepsilon}\mathbb{E}_t\left[\frac{\hat{M}_0}{\hat{S}_0}; \sup_{s\in V_t}|r_t(s)| > \varepsilon\right]$$

$$\leq e^\varepsilon(1+\varepsilon)\mathbb{E}_t\left[\frac{e^{-X}}{S}; X + \log M \geq 0\right] + e^\varepsilon\mathbb{E}_t\left[\frac{M_0}{S_0}; S > (1+\varepsilon)S_0\right]$$

$$+ e^{-\varepsilon}\mathbb{E}_t\left[\frac{\hat{M}_0}{\hat{S}_0}; \sup_{s\in V_t}|r_t(s)| > \varepsilon\right].$$

Therefore,

$$\mathbb{E}_t\left[\frac{e^{-X}}{S}; X + \log M \geq 0\right]$$

$$\geq \frac{e^{-\varepsilon}}{1+\varepsilon}\mathbb{E}_t\left[\frac{e^{-X}}{\hat{S}_0}; X + \log \hat{M}_0 \geq \varepsilon\right]$$

$$- \frac{1}{1+\varepsilon}\mathbb{E}_t\left[\frac{M_0}{S_0}; S > (1+\varepsilon)S_0\right] - \frac{e^{-2\varepsilon}}{1+\varepsilon}\mathbb{E}_t\left[\frac{\hat{M}_0}{\hat{S}_0}; \sup_{s\in V_t}|r_t(s)| > \varepsilon\right].$$

The proof proceeds through a sequence of lemmas, which show that the various error terms are small. Many statements along the way will hold true only up to a constant factor. We use $D, D_1, D_2$, and so on, to denote such (positive) constants. Occasionally the same symbol denotes different constants. □

We start the sequence of lemmas by putting on record the elementary fact that the expectation of a nonnegative random variable over an event can be controlled by the tail of its distribution:

LEMMA 4.2. *For any nonnegative random variable $Y$, any measurable set $A$, and any positive and finite $C$,*

(4.8) $$\mathbb{E}[Y; A] \leq C\mathbb{P}[A] + C\mathbb{P}[Y > C] + \int_C^\infty \mathbb{P}[Y > y]\,dy.$$

PROOF. The proof is elementary and is omitted. □

The second general lemma relates the maxima of a given (deterministic) function to the integral of the same function.

LEMMA 4.3. *Let $h\colon B(r) \to \mathbb{R}$ be a continuously differentiable real-valued function over a closed ball of radius $r$. Let $H = \max_{z\in B(r)}\|\dot{h}(z)\|$. Then, for some positive constant $D$,*

(4.9) $$\int_{B(r)}\exp\left\{h(y) - \max_{z\in B(r)} h(z)\right\}dy \geq (H + D/r)^{-d}.$$

PROOF. The proof is given in the Appendix. □

The next lemma takes us back to the specific terms analyzed in Theorem 4.1:

LEMMA 4.4. *Assume the conditions of Theorem* 4.1. *Then*

$$\mathbb{E}_t\left[\frac{M_0}{S_0}; \sup_{s\in V_t}|r_t(s)| > \varepsilon\right] = o(x^{-\alpha}), \tag{4.10}$$

*and the same holds for* $\hat{M}_0/\hat{S}_0$.

PROOF. The proof is based on Lemma 4.2. Set $A = \{\sup_{s\in V_t}|r_t(s)| > \varepsilon\}$, and $C = n^{D\log n}$, for some positive constant $D$. We begin by showing that

$$n^{D\log n}\mathbb{P}_t\left[\sup_{s\in V_t}|r_t(s)| > \varepsilon\right] = o(x^{-\alpha}). \tag{4.11}$$

Recall that $U_n^{i_1,\ldots,i_k}(s) = \sum_{u\in\mathcal{A}_n}\beta_{u,n}^{i_1,\ldots,i_k}(s)(W_u - \psi'_u(\xi_t\theta_u(t)))$. We write briefly $\beta_{u,n}^{(k)}(\cdot)$ and $U_n^{(k)}(\cdot)$, or even $\beta_{u,n}(\cdot), U_n(\cdot)$ when no confusion is likely. The remainder is absolutely and uniformly bounded in $V_t$ by

$$\begin{aligned}D_1(\log x)^{\alpha+2}x^{-(\alpha+1)}\max_{1\leq i_1,\ldots,i_{\alpha+2}\leq d}\sup_{s\in V_t}|U_n^{(\alpha+2)}(s)|\\ + D_2(\log x)^{\alpha+3}x^{-(\alpha+1)},\end{aligned} \tag{4.12}$$

since $\theta_u(\cdot)$ and its derivatives are uniformly bounded by assumption.

The bound (4.12) on the remainder leads directly to:

$$\mathbb{P}_t\left[\sup_{s\in V_t}|r_t(s)| > \varepsilon\right] \leq \sum\mathbb{P}_t\left[\sup_{s\in V_t}|U_n^{(\alpha+2)}(s)| > Dx(\log x)^{-(\alpha+2)}\right], \tag{4.13}$$

where the sum expands over $\{1,\ldots,d\}^{(\alpha+2)}$. Note that the inequality is valid for any $x$ which is greater than $D(\log x)^{(\alpha+3)}$. The probability on the right hand side of (4.13) is at most the sum of similar tail probabilities of the random field $U_n(s)$ and its negation. These probabilities differ only by a constant. We bound the tail of a bounded field by the expectation of an exponentiated field. Such expectations are investigated next.

Write $F_t(\cdot)$ for the probability distribution function, under the alternative probability $\mathbb{P}_t$, of the collection $\{W_u, u \in \mathcal{A}_n\}$. The expectation $\mathbb{E}_t[\sup_{s\in V_t}\exp\{U_n(s)\}]$ is given, upon dividing and multiplying by $\int_{V_t}\exp\{U_n(s)\}\,ds$ and using Fubini's theorem, by

$$\int_{V_t}\int\left[\frac{\sup_{s\in V_t}e^{U_n(s)}}{\int_{V_t}e^{U_n(s)}\,ds}\right]\times e^{U_n(r)}\,dF_t(w)\,dr, \tag{4.14}$$

where the innermost integral (or a sum for discrete models) is with respect to the sample space. Using a suitable exponential tilting, obtained by considering the cumulant generating function of $\{W_u\}$ under $\mathbb{P}_t$, we have that (4.14)

$$(4.15) \qquad = \int_{V_t} e^{\sum_{u \in \mathcal{A}_n} A_n(r,t)} \times \mathbb{E}_r \left[ \frac{\sup_{s \in V_t} e^{U_n(s)}}{\int_{V_t} e^{U_n(s)} \, ds} \right] dr,$$

for $A_n(r,t) = \psi_u(\xi_t \theta_u(t) + \beta_{u,n}^{(\alpha+2)}(r)) - \psi(\xi_t \theta_u(t)) - \beta_{u,n}^{(\alpha+2)}(r)\psi'_u(\xi_t \theta_u(t))$. Note that the probability distribution $\mathbb{P}_r$, or equivalently its distribution function, belongs as does the probability $\mathbb{P}_t$ to an exponential family. The exact form of the family may be written explicitly but is not essential for the proof.

Taking one additional term in the Taylor expansion for $\psi_u(\xi_t \theta_u(t) + \beta_{u,n}^{(\alpha+2)}(r))$ we can write $A_n(r,t) = \frac{1}{2}(\beta_{u,n}^{(\alpha+2)}(r))^2 \psi''_u(\vartheta)$, for some point $\vartheta$ close to $\xi_t \theta_u(t)$. Then, since $\psi''_u$ are uniformly bounded, we get an upper bound on (4.15):

$$(4.16) \qquad e^D \int_{V_t} \mathbb{E}_r \left[ \frac{\sup_{s \in V_t} e^{U_n(s)}}{\int_{V_t} e^{U_n(s)} \, ds} \right] dr.$$

Define $h(y) = U_n(t + \frac{1}{x} \Lambda_n^{-1/2}(t) y)$ and $B = B(\log x)$. The expectation in (4.16) is, by Lemma 4.3, smaller than or equal to the expectation under $\mathbb{P}_r$ of $(H + D/\log x)^d$, for $H = \max_{z \in B} \| \frac{1}{x} \Lambda_n^{-1/2}(t) \dot{U}_n(t + \frac{1}{x} \Lambda_n^{-1/2}(t) z) \|$. It is sufficient to consider $\mathbb{E}_r(H^d)$, or even $\mathbb{E}_r(\sum_{u \in \mathcal{A}_n} |W_u - \psi'_u(\xi_t \theta_u(t))|)^d$. By independence, and since $W_u^d$ are integrable with respect to $\mathbb{P}_r$, the expectation exhibits a growth which is at most polynomial. Assertion (4.11) now follows from Chebyshev's inequality.

Consider next the tail of $M_0/S_0$. In the application of Lemma 4.3 we may identify $B$ by $B(\log x)$ and $h(y)$ by $xZ_n(t + \frac{1}{x} \Lambda_n^{-1/2}(t) y)$. Then $S_0/M_0$ is in the required form: $|\Lambda_n(t)|^{-1/2} \int_B \exp\{h(y) - \max_{z \in B} h(z)\} \, dy$. It follows that

$$(4.17) \qquad \frac{M_0}{S_0} \leq D(H + D/\log x)^d,$$

where $H \leq D \max_{s \in V_t} \|\dot{Z}_n(s)\| \leq D n^{-1/2} \sum_{u \in \mathcal{A}_n} |W_u|$. By Chebyshev's inequality,

$$(4.18) \qquad \int_{n^{D \log n}}^{\infty} \mathbb{P}_t \left[ \frac{M_0^2}{S_0^2} \geq y^2 \right] dy$$

$$\leq Dn^{-D \log n} \mathbb{E}_t \left( Dn^{-1/2} \sum_{u \in \mathcal{A}_n} |W_u| + D/\log x \right)^{2d}.$$

Clearly, both $C\mathbb{P}_t[M_0/S_0 > C]$, and $\int_C^{\infty} \mathbb{P}_t[M_0/S_0 > y] \, dy$ are $o(x^{-\alpha})$ as requested.

The tail behavior of $\hat{M}_0/\hat{S}_0$ can be managed almost identically. Let $\Xi$ denote the bound (4.12) on the remainder. Hence, $\hat{M}_0/\hat{S}_0 \leq e^{2\Xi} M_0/S_0$, and the problem may be reduced to the evaluation of $\int_{n^{D \log n}}^{\infty} \mathbb{P}_t(e^{2\Xi} \geq y^{1/2}) \, dy$, and $\int_{n^{D \log n}}^{\infty} \mathbb{P}_t(M_0/S_0 \geq y^{1/2}) \, dy$. The latter is clearly in the appropriate order of magnitude [raise each side of inequality (4.17) to the power of 4]. The former can be managed by raising each side to the power of 8, say, and using the method taken in the evaluation of the expectation of the exponentiated field. $\square$

LEMMA 4.5. *Under the conditions of Theorem* 4.1,

$$\text{(4.19)} \qquad \mathbb{E}_t\left[\frac{M_0}{S_0}; S > (1+\varepsilon)S_0\right] = o(x^{-\alpha}).$$

PROOF. The proof's structure is patterned after Lemma 4.2. This time, we set $A = \{S > (1+\varepsilon)S_0\}$ and take as before $C = n^{D \log n}$. Since the radius of a generic ball $V_i$ is $x^{-2}$ we obtain that the total number of balls that are needed in order to cover $T \setminus V_t$ is proportional to $x^{2d}$. However, since the elementary inclusion $\{\sum_{i=1}^n Y_i > y\} \subset \bigcup_{i=1}^n \{Y_i > y p_i\}$ holds for any sequence of nonnegative numbers $p_1, \ldots, p_n$ that add up to 1 and any sequence $Y_1, \ldots, Y_n$ of nonnegative random variables, it is sufficient to evaluate the probability of $\bigcup_{i=1}^{Dx^{2d}}\{S_i > \varepsilon p_i S_0\}$. We take $p_i = x^{-2d}/D$.

A lower bound on $S_0$ may take the form:

$$S_0 \geq D \exp\left\{-D(\log x)^2 x^{-1} \max_{s \in V_t} \|\ddot{U}_n(s)\|\right\}.$$

As for $S_i$, we decompose $X_t(s)$, given the central point $\tau_i$, as:

$$X_t(s) = x(Z_n(s) - Z_n(\tau_i)) + x(Z_n(\tau_i) - Z_n(t)).$$

The first component is bounded by $(\log x) x^{-1} \max_{s \in V_i} \|\dot{U}_n(s)\| + O(1)$. The second component may be decomposed further to produce

$$\text{(4.20)} \qquad x(U_n(\tau_i) - U_n(t)) + x \sum_{u \in \mathcal{A}_n} (\beta_{u,n}(\tau_i) - \beta_{u,n}(t)) \psi_u'(\xi_t \theta_u(t)).$$

The sum in (4.20) is equal to $x(\sum_{u \in \mathcal{A}_n} \beta_{u,n}(\tau_i)\beta_{u,n}(t) - 1)$, up to an $O(x^2 n^{-1/2})$ term. But, $\sum_{u \in \mathcal{A}_n} \beta_{u,n}(\tau_i)\beta_{u,n}(t)$ is (asymptotically) the correlation, under the transformed measure $\mathbb{P}_t$, between $U_n(\tau_i)$ and $U_n(t)$, and as such, it is absolutely bounded by unity. It turns out that the second component in (4.20) is a negative multiple of $x^2$.

Gathering the various random and deterministic terms, one sees that to complete the proof it is sufficient to bound the terms $\mathbb{P}_t(\max_{s \in V_i} \|\dot{U}_n(s)\| \geq x^3/\log x)$, $\mathbb{P}_t(\max_{s \in V_t} \|\ddot{U}_n(s)\| \geq x^3/\log^2 x)$, and $\mathbb{P}_t(U_n(\tau_i) - U_n(t) \geq Dx)$. The probabilities involving the maxima over $V_i$ and $V_t$ can be handled by

the methods used in Lemma 4.4. Notice that the thresholds in this case are even larger than the ones we were using before. The probability associated with $U_n(\tau_i) - U_n(t)$ decays exponentially. This can be verified, for example, by Chebyshev's inequality since the exact form of the cumulant function of $W_u$ is known. $\square$

LEMMA 4.6. *Under the same conditions as before,*

$$\text{(4.21)} \qquad \sum_{i=1}^{K} \mathbb{E}_t\left[\frac{M_i}{S_i}; M_i > 1\right] = o(x^{-\alpha}).$$

PROOF. The handling of $\mathbb{P}_t(M_i > 1)$ may be carried out similarly to Lemma 4.5 by using the decomposition of $X_t(s)$ considered there. The quotient $M_i/S_i$ can be analyzed using Lemma 4.3. $\square$

LEMMA 4.7. *Under the same conditions as before,*

$$\text{(4.22)} \qquad \mathbb{E}_t\left[\frac{e^{-x(Z_n(t)-x)}}{\hat{S}_0}; x(Z_n(t)-x) + \log \hat{M}_0 \geq 0\right] < \infty.$$

PROOF. The investigation of the tail behavior of $\hat{M}_0/\hat{S}_0$ shows that it is integrable with respect to $\mathbb{P}_t$. This is sufficient in order to prove (4.22). $\square$

4.3. *Normal approximation.* Let $\mathcal{U}_n = \mathcal{U}_n(t) = (U_n^{(1)}(t), \ldots, U_n^{(\alpha+1)}(t))^{\text{t}}$ denote the vector of partial derivatives of $U_n(t)$, and set $\mathcal{U}_n = \mathcal{U}_n(t) = (U_n(t), \mathcal{U}_n(t))^{\text{t}}$. We treat $U_n^{(k)}$, the $k$th order partial derivatives, as a $\binom{k+d-1}{k}$-dimensional vector and $\mathcal{U}_n$ is a large column vector. The expectations we would like to evaluate as a result of Theorem 4.1 are of the form $\mathbb{E}_t[G(\mathcal{U}_n)]$, for an appropriate function $G$.

The target at this stage is to show that

$$\text{(4.23)} \qquad \mathbb{E}_t G(\mathcal{U}_n) = \mathbb{E}_t G(\mathcal{U}) P(\mathcal{U}) + o(x^{-\alpha}),$$

for $\mathcal{U} \equiv \mathcal{U}(t) = (U(t), \mathcal{U}(t))^{\text{t}}$, a normally distributed (under $\mathbb{P}_t$) random vector with zero mean and variance–covariance matrix denoted by $C_n$ (see below) and for $P$ an appropriate polynomial. Equation (4.23) corresponds to an application of a central limit theorem to a sum of independent vectors that are not identically distributed. If $\mathcal{U}_n$ possesses a density then one may use Theorem 19.3 of [6] in order to expand the density to the required accuracy and get the result. If $\mathcal{U}_n$ is discrete, as is the case in the numerical example we presented, different tools are required.

We give below an argument corresponding to $\alpha = 1$ and general distributions with a finite forth moment. We shall assume that (a) the functions $\theta_u(\cdot)$, and $\psi_u(\cdot)$ and their derivatives are bounded (see Section 4.2), and that (b) the smallest eigenvalue of the covariance matrix $C_n$ is bounded away from zero. Note that for $\alpha = 1$ we have that $P \equiv 1$.

Let $d_n = D \log n$. The basic strategy we apply in this subsection involves three steps. In the first step, the expectation is evaluated separately on two complementary events, $\{\|\mathcal{U}_n\| \leq d_n\}$ and $\{\|\mathcal{U}_n\| > d_n\}$, with the latter shown to be negligible. The second step utilizes bounds for errors of normal approximation in order to establish (4.23) for expectations restricted to the first event. As a result, $\mathcal{U}_n$ is replaced by $\mathcal{U}$. The third step reverses step one, but in the normal setting. The first step uses techniques similar to those applied in the proof of Lemma 4.4. In [9] one can find a discussion of the expansion of the distribution of the maximum for Gaussian random fields, which involves some modification of the tools that are used in empirical random fields. As a corollary of that discussion, it will be possible to produce a parallel of Lemma 4.4 and to prove the third step. Details for the first and the third steps are omitted. Henceforth, we concentrate on the second step.

We shall need the following notation. Let $b_{u,n}(t) = (\beta_{u,n}(t), \dot{\beta}_{u,n}(t), \ddot{\beta}_{u,n}(t))^t$, where $\ddot{\beta}_{u,n}(t)$ is considered as a vector, and thus $\mathcal{U}_n(t) = \sum_u b_{u,n}(t)[W_u - \psi'_u]$, for $\psi'_u = \psi'_u(\xi_t \theta_u(t))$. Define $C_n = \sum_u \mathrm{Cov}(b_{u,n}(t)[W_u - \psi'_u])$, and let

$$X_u = X_{u,n}(t) = n^{1/2} C_n^{-1/2} b_{u,n}(t)[W_u - \psi'_u].$$

Let $\mathcal{V}_n(t) = \sum_u X_u = n^{1/2} C_n^{-1/2} \mathcal{U}_n(t)$ be the sum of the independent vectors and let $Q_n$ be the distribution of $\mathcal{W}_n(t) = n^{-1/2} \mathcal{V}_n(t) = C_n^{-1/2} \mathcal{U}_n(t)$.

With those notations in mind, write

$$\begin{aligned}(4.24) \quad & \mathbb{E}_t[G(\mathcal{U}_n); \|\mathcal{U}_n\| \leq d_n] - \mathbb{E}_t[G(\mathcal{U}); \|\mathcal{U}\| \leq d_n] \\ & = \int f(w) \, d(Q_n(w) - \Phi(w)),\end{aligned}$$

where $\Phi$ is the Gaussian measure, and $f(w) = G(C_n^{1/2} w) I_{\{\|C_n^{1/2} w\| \leq d_n\}}$.

The term in (4.24) may be bounded with the aid of Theorem 13.3 of [6], which states that

$$(4.25) \quad \left| \int f \, d(Q_n - \Phi) \right| \leq \omega_f(\mathbb{R}^k) a_3(k) \rho_4 n^{-1/2} \\ + \tfrac{4}{3} \omega_f^*(2^{7/2} \pi^{-1/3} k^{4/3} \rho_3 n^{-1/2}; \Phi),$$

where $a_3(k)$ is a constant that depends only on $k$, the dimension of $\mathcal{U}_n$ and, for $\rho_i = n^{-1} \sum_u \mathbb{E}\|X_u\|^i$, $i = 3, 4$. The term $\omega_f(A) = \sup\{|f(x) - f(y)| : x, y \in A\}$ is the modulus of oscillation of the function $f$ over $A$ and, for $\varepsilon > 0$ and a measure $\mu$,

$$\omega_f^*(\varepsilon; \mu) = \sup_{y \in \mathbb{R}^k} \int \omega_f(B(\varepsilon) + x - y) \mu(dx).$$

The fourth moment is bounded by $n\|C_n^{-1/2}\|^4 \sum_u \|b_{u,n}(t)\|^4 \mathbb{E}_t[W_u - \psi'_u]^4$. Under the assumption regarding the boundedness of $\theta_u(\cdot)$, $\psi_u(\cdot)$, and their derivatives, we have that $n\sum_u \|b_{u,n}(t)\|^4 \mathbb{E}_t[W_u - \psi'_u]^4 = \Theta(1)$. Moreover, $\|C_n^{-1/2}\|^4 = (\lambda_{\min}(C_n))^{-2}$, where $\lambda_{\min}(C_n)$ is the smallest eigenvalue of $C_n$, which is bounded away from zero. Thus, the (averaged) fourth moment is finite.

In the assessment of $\omega_f$ and $\omega_f^*$ it is enough to consider $f(\mathfrak{u})$ for $\mathfrak{u} = C_n^{-1/2} w$, since $C_n$ is an $\Theta(1)$ matrix. For such $\mathfrak{u}$,

$$f(\mathfrak{u}) \leq e^\varepsilon \frac{\hat{M}_0}{\hat{S}_0} I_{\{\|\mathfrak{u}\| \leq d_n\}}.$$

We regard here $\mathcal{U}_n$ as fixed at the value $\mathfrak{u} = (u, \dot{u}, \ddot{u})$. Direct maximization, and the use of elementary inequalities [such as $F_{\chi^2(d,\theta^2)}(y^2) \geq e^{-\theta^2/2} F_{\chi^2(d)}(y^2)$, which features a relationship between the c.d.f. of a noncentral chi-squared distribution and a central one] lead to

$$(4.26) \quad \frac{\hat{M}_0}{\hat{S}_0} I_{\{\|\mathfrak{u}\| \leq d_n\}} \leq D e^{\dot{u}' \Lambda_n^{-1}(t) \dot{u}/2} I_{\{\|u\| \leq d_n\}} \leq D n^{1/4} I_{\{\|\mathfrak{u}\| \leq d_n\}}.$$

Therefore, $f$ is bounded, and the first term in the right-hand side of (4.25) is of the order of magnitude of $O(n^{-1/4})$.

Consider next $w_f^*$. Let $A$ be any subset of $\mathbb{R}^k$. The modulus of oscillation of the indicator function, $I_A$, over the ball $B(\varepsilon) + x$, is

$$(4.27) \quad w_{I_A}^*(\varepsilon; \Phi) = \sup_{y \in \mathbb{R}^k} \int I_{(\partial A)^\varepsilon}(x - y) \Phi(dx) = \sup_{y \in \mathbb{R}^k} \Phi((\partial(A + y))^\varepsilon),$$

where $A^\varepsilon$ is the set of points whose distance from $A$ is less than $\varepsilon$. Applying Corollary 3.2 of [6] (with $s = 0$ and $\rho = \varepsilon$) we obtain that for a convex set $A$, $w_{I_A}^*(\varepsilon; \Phi) \leq b(k)\varepsilon$, where $b(k)$ is a constant depending on $k$ only. Hence, since the set $\{\|u\| \leq d_n\}$ is convex and $f$ is bounded, we conclude that

$$w_f^*(\varepsilon; \Phi) \leq D n^{1/4} b(k) \varepsilon = O(n^{-1/4}).$$

This completes the proof that (4.25) is bounded by $O(n^{-1/4}) = o(1/x)$.

4.4. *Elimination of the indicator.* This section proceeds the analysis of the two (compactly written) terms,

$$\mathbb{E}_t[e^{-X}/\hat{S}_0; X + \log \hat{M}_0 \geq \pm \varepsilon],$$

which originated from the analysis in the previous sections. In order to remove the dependency on the event, we proceed in two steps. The first step involves a conditioning argument, where the conditioning here is with respect to the $\sigma$-field generated by the derivative components $\dot{\mathcal{U}}$. These components generate the local random field. Note that both $\hat{S}_0$, and $\hat{M}_0$ are measurable with respect to this $\sigma$-field. The second step presents a Mill's ratio type of an approximation of the innermost expectation. Recall that $(U, \mathcal{U})$ are jointly Gaussian, hence the conditional distribution of $U$, given the $\sigma$-field, is normal. We apply a simple lemma which is proved in the Appendix:

LEMMA 4.8. *Let $Y \sim N(\mu, \sigma^2)$. Suppose $x \to \infty$ and $y \to 0$, such that $xy$ converges to a constant. Then,*

$$\mathbb{E}[e^{-xY}; Y \geq y]$$
$$= e^{-xy} \cdot \phi\left(\frac{y-\mu}{\sigma}\right) \cdot \sum_{m=0}^{k} \frac{(-1)^m (2m)!}{\sigma^{2m} 2^m m!} \left[x + \frac{y-\mu}{\sigma^2}\right]^{-(2m+1)} + r,$$

*where $\phi(\cdot)$ is the standard normal density function and $r = o([x + \frac{y-\mu}{\sigma^2}]^{-(2k+2)})$.*

Let $\mathcal{F}_t = \sigma(\mathcal{U})$ denote the $\sigma$-field generated by the $k$th order partial derivatives of the (centered) field at $t$, for $1 \leq k \leq \alpha + 1$ and for some $\alpha \leq 2$. Let $Z(t) = U(t) + x$ be a random variable which, under $\mathbb{P}_t$, is normally distributed with mean $x$ and variance $\text{Var}_t(U(t))$. Let $\rho = \text{Cov}(Z(t), \mathcal{U})$ and $\Sigma = \text{Cov}(\mathcal{U})$ [recall that $C_n = \text{Cov}(\mathcal{U})$]. Observe that the conditional mean of $Z(t) - x$ is

(4.28) $$\mu = \mathbb{E}_t(Z(t) - x | \mathcal{F}_t) = \langle \rho, \Sigma^{-1} \mathcal{U} \rangle,$$

while the conditional variance is $\sigma^2 = \text{Var}_t(Z(t) - x | \mathcal{F}_t) = 1 - \langle \rho, \Sigma^{-1} \rho \rangle$. Therefore, conditioning on $\mathcal{F}_t$, taking $y = (-\log \hat{M}_0 \pm \varepsilon)/x - r_n(t)$ and applying Lemma 4.8 we obtain, up to a $o(x^{-\alpha})$ remainder,

$$\mathbb{E}_t\left[\frac{e^{-X}}{\hat{S}_0}; X + \log \hat{M}_0 \geq \pm \varepsilon\right]$$

(4.29) $$= \frac{e^{\mp \varepsilon}}{x} \mathbb{E}_t\left[\frac{\hat{M}_0}{\hat{S}_0} \cdot \frac{e^{-1/(2\sigma^2)(\mu + (\log \hat{M}_0/x) \mp \varepsilon/x + r_n(t))^2}}{\sqrt{2\pi\sigma^2}}\right.$$
$$\left. \times \sum_{m=0}^{\lfloor \alpha/2 \rfloor} \frac{[(-1)^m (2m)!]/[x^{2m} \sigma^{2m} 2^m m!]}{(1 - [\mu/x + (\log \hat{M}_0/x^2) \mp \varepsilon/x^2 + r_n(t)/x]/\sigma^2)^{2m+1}}\right].$$

The difference $r_n(t) = \sum_{u \in \mathcal{A}_n} \beta_{u,n}(t) \psi'_u(\xi_t \theta_u(t)) - x$ was defined in (1.11).

Although the expression in square brackets is pretty involved, note that the dependency on the event is removed and the result is an expectation of a function of the derivatives. The next section shows a general algorithm to produce an explicit approximate evaluation of this expectation.

4.5. *Evaluation of the functional.* The expression in the square brackets on the right-hand side of approximation (4.29) has four significant terms: $\hat{M}_0$, $\hat{S}_0$, $\phi(-\log \hat{M}_0/x\sigma - \mu/\sigma - r_n(t)/\sigma)$, and the finite sum of rational functions. In order to obtain an approximation with error term of the order of magnitude of $o(x^{-\alpha})$ we show below how each term can be expanded as power series in $1/x$. Clearly, $\varepsilon$ can be ignored in what follows. Except for the second term, $\hat{S}_0$, the three other involve the maxima $\hat{M}_0$. We leave the discussion about $\hat{S}_0$ to the end of the section, and commence with the maxima.

For each $t$, the function $\log \hat{M}_0$ is the maximum, subject to constraints, of a polynomial. The coefficients are functions of the random and deterministic derivatives. The constraint is the neighborhood $V_t$ of $t$. The dominant part in the polynomial $\hat{X}_t(s)$ are the linear component associated with the random gradient and the quadratic component associated with the deterministic Hessian. We ignore here the issue of the constraint, since it does not affect the evaluation to the order of accuracy considered here, and proceed with the unconstrained maximization. We reintroduce maximization under constraints when we deal with boundary effects, for which constraints are significant and may be active.

The change of variable $y = x(s-t)$ will produce the representation

$$f(y) = \hat{X}_t(y/x + t)$$
$$= \langle y, \dot{U}(t) \rangle - \tfrac{1}{2} \langle y, \Lambda_n(t) y \rangle + \sum_{j=1}^{\alpha} x^{-j}[g_{j+1}(y, U) + h_{j+2}(y, \beta)],$$

where

$$g_j(y, U) = y_{i_1,\ldots,i_j} U^{i_1,\ldots,i_j}(t)/j!,$$
$$h_j(y, \beta) = y_{i_1,\ldots,i_j} \sum_{u \in \mathcal{A}_n} \beta_{u,n}^{i_1,\ldots,i_j} \psi'_u(\xi_t \theta_u(t))/xj!.$$

The remainder terms $r_1, r_2$ that appear in (4.6) are discarded in the above but may be introduced. In general, the level of accuracy $\alpha$ determines whether remainder terms (which depend both on $x$ and on $n$) will affect the final result or not. Let $f(y) = q(y) + r(y)$, where $q(y)$ is the dominant quadratic part, and $r(y)$ is the remainder. We denote the maximum value $\log \hat{M}_0$ by $f(y^*)$, for $y^*$ the maximizer. In order to approximate this maximum we use a quasi Newton–Raphson algorithm (see, e.g., [8]). This entails finding a point $\hat{y}$ which approximates the location of the maximum of $f(y)$ up to a certain order. This point, in turn, will also be used to bound the difference $f(\hat{y}) - f(y^*)$.

To be more specific, take $y_0 = 0$, $y_1 = \Lambda_n^{-1}(t)\dot{U}(t)$, and define recursively

$$y_{k+1} = y_k + \Lambda_n^{-1}(t)\dot{f}(y_k) = \Lambda_n^{-1}(t)\dot{U}(t) + \Lambda_n^{-1}(t)\dot{r}(y_k).$$

Plugging this back into $f$ produces

$$f(y_{k+1}) = \tfrac{1}{2}\langle \dot{U}(t), \Lambda_n^{-1}(t)\dot{U}(t)\rangle$$
$$- \tfrac{1}{2}\langle \dot{r}(y_k), \Lambda_n^{-1}(t)\dot{r}(y_k)\rangle + r(\Lambda_n^{-1}(t)[\dot{U}(t) + \dot{r}(y_k)]).$$

The incremental increase is

$$f(y_{k+1}) - f(y_k)$$
$$= \langle \dot{r}(y_k) - \dot{r}(y_{k-1}), \Lambda_n^{-1}(t)[\dot{r}(\tilde{y}) - (\dot{r}(y_{k-1}) + \dot{r}(y_k))/2]\rangle,$$

for some $\tilde{y}$ on the line segment connecting $y_k$ and $y_{k-1}$.

By running the iterative process $\lfloor \alpha/2 \rfloor + 1$ times, one gets that the difference $f(\hat{y}) - f(y^*)$ is $o(x^{-\alpha})$, as required. By the form of the function $f$, and the definition of the recursive Newton–Raphson sequence, it is clear that $f(\hat{y})$ is also a polynomial, which can be rearranged according to powers of $1/x$. The coefficient of the polynomial are sums of products of both the random and deterministic derivatives.

The finite sum, which appears on the right-hand side of expression (4.29), is clearly a functional of the derivatives only; it contains $\mu$, which is a linear function of the derivatives vector $\mathcal{U}$, and $\log \hat{M}_0$. After substituting $f(\hat{y})$ for $\log \hat{M}_0$, this sum is expressible as power series. Terms with smaller order of magnitude than $x^{-\alpha}$ can be ignored. This way, a finite representation in the form of a polynomial may be obtained.

Another term involves the normal density. A Taylor expansion of the exponent function, combined with the polynomial representation of $\log \hat{M}_0$, enables us to rewrite the term as $\frac{1}{\sigma}\phi(\mu/\sigma)$ multiplied by a polynomial, with coefficient which are functions of the derivatives. In the final evaluation of the expectation, the normal density $\frac{1}{\sigma}\phi(\mu/\sigma)$ may be absorbed back into the joint normal density of the partial derivatives vector $\mathcal{U}$. This leads to another normal density function, easily recognized as the conditional density of $\mathcal{U}$ given $\{Z(t) = x\}$. To see this, note that $\mathbb{E}_t(\mathcal{U}|\{Z(t) = x\}) = 0$, $\mathrm{Var}_t(\mathcal{U}|\{Z(t) = x\}) = \Sigma - \rho \otimes \rho$, and that the following two relations hold:

$$\frac{\mu^2}{\sigma^2} + \langle \mathcal{U}, \Sigma^{-1}\mathcal{U}\rangle = \left\langle \mathcal{U}, \left[\Sigma^{-1} + \frac{\Sigma^{-1}[\rho \otimes \rho]\Sigma^{-1}}{1 - \langle \rho, \Sigma^{-1}\rho\rangle}\right]\mathcal{U}\right\rangle$$
$$= \langle \mathcal{U}, [(\Sigma - \rho \otimes \rho)^{-1}]\mathcal{U}\rangle$$

and

$$\det(\Sigma - \rho \otimes \rho) = \det(\Sigma) \cdot (1 - \langle \rho, \Sigma^{-1}\rho\rangle) = \sigma^2 \det(\Sigma).$$

Finally, we consider the ratio between $\hat{M}_0$ and $\hat{S}_0$, which is given by the integral of $\exp\{f(y) - f(y^*) + \delta(t + y/x) - \delta(t)\}$ over the region $\Lambda_n^{-1/2}(t)B(\log x)$.

By centering about $y_1$, the first point in the Newton–Raphson series, we obtain the following approximation:

$$(2\pi)^{d/2}|\Lambda_n(t)|^{-1/2}$$
$$\times \left\{ 1 + \mathbb{E}_Y[r(Y) - r(y_1) + \delta(t + Y/x) - \delta(t)] \right.$$
$$\left. + \cdots + \frac{1}{\alpha!}\mathbb{E}_Y[r(Y) - r(y_1) + \delta(t + Y/x) - \delta(t)]^\alpha \right\},$$

for $Y \sim N(y_1, \Lambda_n^{-1}(t))$. The integrands, $(r(Y) - r(y_1))^m$, are polynomial in $1/x$, with coefficients that are polynomials of the partial derivatives. The expectation of moments of polynomials will produce again polynomials. The term $\delta(t + y/x) - \delta(t)$ should be expanded about $t$ up to the required order of accuracy. Specifically, if we let $x = (n^{1/\nu})$, for $4 \leq \nu < 6$, then, conditional on $Y = y$, we have the expansion

$$(4.30) \qquad \frac{\langle y, \dot{\delta}(t) \rangle}{x} + \frac{1}{2}\frac{\langle y, \ddot{\delta}(t)y \rangle}{x^2} + \cdots + \frac{1}{m!}\frac{y^{i_1,\ldots,i_m}\delta_{i_1,\ldots,i_m}(t)}{x^m},$$

where $m = \lceil 3 + \alpha - \nu/2 \rceil - 1$. For $\nu \geq 6$ (and any given accuracy level $\alpha$), $\delta(t + y/x) - \delta(t)$ may be approximated by zero since $\delta(t)$ is of the order of magnitude of $O(x^3 n^{-1/2})$.

Incorporating all the approximating expressions for the four terms mentioned earlier (by standard polynomial multiplication and division), one obtains an $\alpha$-degree polynomial in $1/x$. The coefficients of the polynomial are of the form of a sum of products of the partial derivatives. The expectation of these products, with respect to the conditional distribution of $\dot{\mathcal{U}}$, can be handled using Wick's formula (see Adler [2]).

4.5.1. *Boundary effect.* Up to this point the derivation of the terms associated with the local behavior of the random field ignored the boundedness of the set of parameters $T$. This is justified at points $t$ for which the local neighborhood $V_t$ is a subset of $T$. However, the approximation as presented above does not apply at points near the boundary, where $V_t$ extends beyond $T$. In this section we outline the modifications that allow a rigorous expansion in the vicinity of the boundary. Indeed, substituting $V_t \cap T$ for $V_t$ and walking once more the path that brought us to this point, one can observe that the proofs hold, essentially word for word, up to Section 4.5. No more than simple regularity conditions regarding the smoothness of the boundary are needed. Divergence in the details of the argument occur in Section 4.5, which deals with the presentation of the relevant term in the form of a conditional expectation of a function of the derivatives of the local process. Again, this function is composed of the elements that resulted from the Mill's ratio type of approximation and the ratio between the maximum of the exponentiated local field and the integral of that field. However, maximization and integration now occur only over the constrained region $V_t \cap T$.

Recall that the parameter set was represented in (1.7) in the form $T = \{t \in \mathbb{R}^d : g_i(t) \leq 0, 1 \leq i \leq m\}$, for a finite collection of smooth constraint functions. Here we develop the algorithms for obtaining approximations of the constrained maximum and integral for this representation. In particular, in this section we will assume that $t = t(x)$ is such that $x \cdot g_j(t)$ converges to a negative constant for a given $j$, whereas $x \cdot g_i(t) \leq -\log x$, for $i \neq j$. Depending on the order of approximation required, other situations may need to be considered. The algorithms presented for the present context can be generalized in a straightforward way in order to deal with other situations as well. Some care must be taken with the difference $\delta_n(t + y/x) - \delta_n(t)$ and for the remainder $r_n(t)$.

Let us start with the maximum of the local field. Using the notation of Section 4.5 we will consider a function of the form $f(y) = q(y) + r(y)$, where $q(y)$ is the quadratic dominant part, and $r(y)$ is the remainder. Dropping the index $j$, we represent the constraint function in the form

$$xg(t + y/x) = x \cdot g(t) + \langle \dot{g}(t), y \rangle + x \cdot v(t, y) = g + \langle \dot{g}, y \rangle + v(y),$$

for an appropriate remainder function $v$. Note, that the linear part is dominant over the remainder part. Since the target function $f$ is asymptotically concave, it follows that the constraint is active, unless the global maximum is in the interior of $V_t \cap T$. The global maximum $\hat{M}_0$ was approximated in the previous section. What is remains is to compute the maximum when the constraint function is active. The Sequential Quadratic Programming (SQP) algorithm (see [8]) is a natural generalization of the quasi Newton–Raphson algorithm from the previous section to the case where the constraint is active. The Lagrangian is maximized by an iterative refinement of its pair of arguments $(y, \lambda)$. The SQP algorithm applies the recursive formula:

$$\begin{pmatrix} y_{k+1} \\ \lambda_{k+1} \end{pmatrix} = \begin{pmatrix} y_k \\ \lambda_k \end{pmatrix} - H^{-1} \begin{pmatrix} \dot{r}(y_k) + \lambda_k \dot{v}(y_k) \\ v(y_k) \end{pmatrix},$$

where

$$-H^{-1} = \begin{pmatrix} \Lambda_n^{-1}(t) & 0 \\ 0 & 0 \end{pmatrix}$$
$$- \frac{1}{\langle \dot{g}, \Lambda_n^{-1}(t)\dot{g} \rangle} \begin{pmatrix} \Lambda_n^{-1}(t)[\dot{g} \otimes \dot{g}']\Lambda_n^{-1}(t) & \Lambda_n^{-1}(t)\dot{g} \\ \dot{g}^{\mathrm{t}}\Lambda_n^{-1}(t) & 1 \end{pmatrix},$$

Starting at the $(d+1)$-origin $(0,0)$, under appropriate regularity conditions, one can show that $f(y_k)$ approximates the constraint maximum up to an order of magnitude of $1/x^k$.

The term $\hat{S}_0$ was approximated in the previous section by the expectation of polynomials of a Gaussian random vector $Y$ with mean vector $y_1$ and covariance matrix $\Lambda_n^{-1}(t)$. In the presence of a boundary, the expectation should be replaced by the expectation over the event $\{g + \dot{g}'Y + v(Y) \leq 0\}$. Below we indicate how one can approximate such expectations.

Observe that the vector $Y$ may be decomposed into two orthogonal (hence independent) components: $\dot{g}'Y = Y_1 \sim N(\dot{g}'y_1, \langle \dot{g}, \Lambda_n^{-1}(t)\dot{g}\rangle)$ and $[I - \Lambda_n^{-1}(t)[\dot{g} \otimes \dot{g}]/\langle \dot{g}, \Lambda_n^{-1}(t)\dot{g}\rangle]Y = Y_2$. The polynomials can then be reformulated in terms of polynomials in $Y_1$ and $Y_2$ and the event can be represented in the form: $\{Y_1 \leq -g - v(Y_1, Y_2)\}$. This event can be approximated by the event $\{Y_1 \leq \tilde{v}(g, Y_2)\}$, for $\tilde{v}$ formed by collection of all terms of the appropriate order from the iterative application of the function $G(z_1) = -z_1 - v(z_1, Y_2)$, starting at the point $z_1 = g$. The following step is the computation of the expectation with respect to $Y_1$. One can use the recursion $\psi_j(y) = -y^{j-1}\phi(y) + (j-1)\psi_{j-2}(y)$, where $\psi_j(y) = \int_{-\infty}^y z^j \phi(z)\,dz$, $\phi$ the standard normal density. The recursion is initiated by $\psi_0(y) = \Phi(y)$, the normal c.d.f. function, and $\psi_1(y) = -\phi(y)$. Finally, after a Taylor expansion of the outcome, an expectation is taken with respect to $Y_2$. This involves expectation of a Gaussian polynomial, and may be carried out with the aid of Wick's formula.

4.5.2. *A detailed example.* Let us consider the exact form of the expansion for $\alpha = 1$, and $x = o(n^{1/4})$. Putting together the previous arguments we have a second order expansion of the form:

$$\tag{4.31} x^{d-1}(2\pi)^{-d/2}\phi(x)\int_T e^{-\delta_n(t)}|\Lambda_n(t)|^{1/2} \\ \times \left[1 - r_n^2(t)/2\sigma_n^2(t) + \frac{1}{x}\tilde{\mathbb{E}}_t(a_1) + \tilde{\mathbb{E}}_t(\bar{a}_1)\right]dt.$$

The expectation $\tilde{\mathbb{E}}_t$ here is with respect to the conditional distribution of $\mathcal{U}$ given $\{Z(t) = x\}$. The coefficient $a_1$ is given by

$$a_1 = \mathbb{E}_Y\langle Y, \dot{\delta}_n(t)\rangle + \langle \rho, \Sigma^{-1}\mathcal{U}\rangle(1 - \|\Lambda_n^{-1/2}(t)\dot{U}(t)\|^2/2)/\sigma_n^2(t) \\ - [(\mathbb{E}_Y g_2(Y, \mathcal{U}) - g_2(y_1, \mathcal{U})) + (\mathbb{E}_Y h_3(Y, \beta) - h_3(y_1, \beta))].$$

It can be shown that all terms in the definition of $a_1$ but the first one are centered variables. Moreover, an expectation of $\langle Y, \dot{\delta}_n(t)\rangle$ with respect to $\mathbb{P}_Y$, followed by an expectation with respect to $\tilde{\mathbb{P}}_t$ shows that this term vanishes as well, so $\tilde{\mathbb{E}}_t(a_1) = 0$.

For the boundary we have a parallel result with regard to $\bar{a}_1$:

$$\bar{a}_1 = \left[I_A + I_{A^c} \cdot \exp\left\{-\frac{1}{2}\left(\frac{y\|\dot{g}_t\| - \langle \dot{g}_t, y_1\rangle}{(\langle \dot{g}(t), \Lambda_n^{-1}(t)\dot{g}(t)\rangle)^{1/2}}\right)^2\right\}\right] \\ \times \left[\Phi\left(\frac{y\|\dot{g}(t)\| - \langle \dot{g}(t), y_1\rangle}{(\langle \dot{g}(t), \Lambda_n^{-1}(t)\dot{g}(t)\rangle)^{1/2}}\right) - 1\right]^{-1},$$

where $A = \{\langle \dot{g}(t), y_1 \rangle \leq y \|\dot{g}(t)\|\}$. Integration along the boundary of $T$ produces, after some manipulations,

$$\frac{\phi(x)x^{d-2}}{(2\pi)^{d/2}} \int_{\partial T} \int_0^\infty e^{-\delta_n(t)} |\Lambda_n(t)|^{1/2} \tilde{\mathbb{E}}_t(\bar{a}_1) \, dy \, dV_{\partial T}(t)$$

$$= \frac{\phi(x)x^{d-2}}{2(2\pi)^{(d-1)/2}}$$

$$\times \int_{\partial T} [e^{-\delta_n(t)} |\Lambda_n(t)|^{1/2} \langle \dot{g}(t), \Lambda_n^{-1}(t)\dot{g}(t) \rangle^{1/2} / \|\dot{g}(t)\|] \, dV_{\partial T}(t).$$

To summarize, we obtain:

THEOREM 4.9. *Let $\{Z_n(t); t \in T\}$ be a random field given, marginally, by* (1.4). *Assume that $T$ is a compact, convex subset of the $d$-dimensional Euclidean space. Let the conditions of Theorem* 4.1 *with $\alpha = 1$ hold. Assume further that $\theta_u(\cdot)$, and $\psi_u(\cdot)$ are bounded, as well as their second and forth order derivatives, respectively. Finally, assume that the function $g$, which defines the boundary of $T$, is piecewise continuously differentiable. Let $x \to \infty, n \to \infty$ be such that $x = o(n^{1/4})$. Then $\mathbb{P}(\sup_{t \in T} Z_n(t) \geq x)$ is approximated by* (1.6). *For $x = o(n^{1/5})$ the term $r_n^2(t)/2\sigma_n^2(t)$ can be neglected.*

**5. Extensions.** In this work we presented a second order approximation, which accounts for edge effects, for a specific class of random fields. At most locations it was indicated in the proofs how one may carry out a higher order approximation by including more terms in the Taylor expansion of the local field and by using higher order asymptotic expansions (for continuous distributions). The product is a functional of these derivatives under a Gaussian joint distribution, which needs to be evaluated. Some details, such as the appropriate extension of the Mill's ratio expansion, have been omitted. Additional work would also be required in order to obtain higher order boundary corrections, which take into account curvature of the boundary and points of nondifferentiability. Although higher order expansions may produce better approximations, it is likely that they will be very complex and not provide additional insight. In addition, numerical examples in the Gaussian case suggest that the two term approximation is frequently reasonably accurate.

The methods we have developed are quite general. The "moving average" representation we have assumed for $X_t$ has much the same effect that moving average representations have long assumed in time series analysis. It allows us to derive by calculation detailed estimates needed in our arguments, especially our use of versions of a local central limit theorem to prove (4.23). It is easy to identify places in our argument, where one could rewrite the result of a calculation as an assumption to provide an approximation for a more general class of random fields, which then could be specialized to subclasses by checking the assumptions.

For the well developed case of a Gaussian random field $\{Z(t) : t \in T\}$, with components that are standard normal, several technical conditions and aspects of the proof can be significantly improved and simplified. See [9] for details.

An interesting generalization is random fields, even Gaussian random fields, with both smooth and nonsmooth components. Such fields arise naturally in the monitoring of images over time, so time is another component that should be added to the parameter space. At each given point in time the score may vary smoothly as the function of the structure of the signal. However, if a signal may abruptly appear, then the score will not be smooth in the direction of the time component. The advantage of our method, which was originally developed in nonsmooth settings, is its flexibility. The essential argument is blind with respect to issues of smoothness, and the calculations can be carried out in a unified manner. Only in the detailed investigation of the local field does the level of smoothness become important.

## APPENDIX: MORE PROOFS

PROOF OF LEMMA 4.3. Let $B = B(r)$ and assume that $\hat{y} \in [-r, r]$ maximizes $h$. The one-dimensional and multi-dimensional cases are dealt separately. Consider first the case $d = 1$. If $\hat{y} = r$ then $h(y) - h(\hat{y}) \geq H(y - \hat{y})$, which gives

(A.1)
$$\int_{-r}^{r} e^{h(y)-h(\hat{y})} \, dy \geq \int_{0}^{2r} e^{-Hy} \, dy$$
$$= H^{-1}(1 - e^{-2rH}) \geq (H + 1/2r)^{-1}.$$

The last inequality is verified by recalling the elementary inequality $y/(1 - e^{-y}) \leq y + 1$, which holds for every $y > 0$. The other extremal situation, $\hat{y} = -r$, is exactly the same. One only has to consider the steepest negative slope $-H$ instead of $H$. If $|\hat{y}| < r$ we use them both to bound $h(y) - h(\hat{y})$ from below by $-H|y - \hat{y}|$. The integral is then evaluated over $[-r, \hat{y}]$ for positive $\hat{y}$, and over $[\hat{y}, r]$ for negative $\hat{y}$. Consequently, a lower bound for interior points, $(H + 1/r)^{-1}$, is achieved, which together with the right-hand side of (A.1) completes the proof for the one-dimensional case.

The multi-dimensional case is again divided into two scenarios. One is when $\hat{y}$ lies on the boundary and the other is when $\hat{y}$ is an interior point. We first treat boundary points $\|\hat{y}\| = r$. To begin with, assume that the projections of $\hat{y}$ onto each plane $y_i - y_j$, $i \neq j$, is the bisector of the axes $y_i, y_j$ and further, that $\hat{y}_i > 0$, for $1 \leq i \leq d$. That is, assume that the $\hat{y}_i$'s are positive and that $\hat{y}_1 = \hat{y}_2 = \cdots = \hat{y}_d$, with a common value which must be equal to $rd^{-1/2}$. We refer to this $\hat{y}$ as "simple." For every $y \in B$ there exists, by the mean-value theorem, a point $\eta$ on the line segment joining $y$ and $\hat{y}$, such that $h(y) - h(\hat{y}) = (y - \hat{y})'\dot{h}(\eta)$. Now, denote by $Q = \prod_{i=1}^{d} Q_i$ the ($d$-dimensional) cube inscribed in $B$ having sides $2rd^{-1/2}$. Therefore,

$$(A.2) \qquad \int_B e^{h(y)-h(\hat{y})}\,dy \geq \prod_{i=1}^{d} \int_{Q_i} e^{(y_i-\hat{y}_i)\dot{h}_i(\eta)}\,dy_i.$$

By a reduction to the one-dimensional case it is clear that every univariate integral above is in the form considered there for boundary points, with $r$ replaced by $rd^{-1/2}$. Since $H_i \equiv \max_{z \in B} |\dot{h}_i(z)| \leq H$ for all $1 \leq i \leq d$, we obtain a lower bound $(H + 1/2rd^{-1/2})^{-d}$. The general case of a boundary point $\hat{y}$, which is not in the form depicted above, is handled similarly except that the axes should be rotated first, to make $\hat{y}$ "simple." Denote by $\Gamma$ the corresponding rotation matrix. Since $\Gamma$ is orthogonal norm preserving transformation the bound changes only slightly, and is given in the form indicated by the lemma.

Last, let $\hat{y}$ be an interior point. Assume, with no loss of generality, that $\hat{y}$ is simple (and also $\hat{y}_i < rd^{-1/2}$). This is clearly enough by the above argument. We repeat the scheme that yielded relation (A.2). This time each univariate integral is evaluated over the subinterval $[-rd^{-1/2}, \hat{y}_i]$. This is akin to the one-dimensional case with $rd^{-1/2}$ substituted for $r$, and thus, the requested bound is obtained. □

PROOF OF LEMMA 4.8. Consider first the case $\sigma^2 = 1$. The general case will follow by substituting $\sigma x$ for $x$, $\mu/\sigma$ for $\mu$, and $y/\sigma$ for $y$. The log-likelihood ratio associated with transforming from mean $\mu$ to mean $y$ is $(\mu - y)Y - (\mu^2 - y^2)/2$. Therefore,

$$\mathbb{E}_\mu(e^{-xY}; Y \geq y) = \mathbb{E}_y(e^{-(x+y-\mu)(Y-y)}; Y - y \geq 0) \times e^{-xy-(y-\mu)^2/2}.$$

The random variable $Y - y$ has a standard normal distribution under the distribution measure $\mathbb{P}_y$. Its density satisfies:

$$\left| \phi(z) - \frac{1}{\sqrt{2\pi}} \sum_{m=0}^{k} \frac{(-1)^m z^{2m}}{2^m m!} \right| \leq \frac{z^{2n+2}}{2^{n+1}(n+1)!}.$$

The approximation and bound on the error are obtained by integration. □

**Acknowledgments.** The authors thank two referees for very helpful comments.

Y. Nardi  
Department of Statistics  
Carnegie Mellon University  
Pittsburgh, Pennsylvania  
USA  
E-mail: yuval@stat.cmu.edu

D. O. Siegmund  
Department of Statistics  
Stanford University  
Stanford, California  
USA  
E-mail: siegmund@stat.stanford.edu

B. Yakir  
Department of Statistics  
Hebrew University  
Jerusalem  
Israel  
E-mail: msby@mscc.huji.ac.il